\documentclass[12pt]{article}
\usepackage{amsmath,amssymb}
\newtheorem{thm}{Theorem}[section]
\newtheorem{lem}[thm]{Lemma}

\newcommand{\R}{{\mathbb R}}

\newcommand{\old}[1]{}
\newcommand{\Adj}{{\rm Adj}}
\newcommand{\C}{{\mathbb C}}
\newcommand{\Z}{{\mathbb Z}}
\newcommand{\M}{{\cal M}}

\newcommand{\Gd}{{\cal G}_D}
\newcommand{\Gdr}{{\cal G}_D^r}

\newcommand{\Gt}{{\cal G}_T}
\newcommand{\Tr}{\mathcal T}
\newcommand{\F} {\mathcal F}
\newcommand{\1} {\overline{1}}
\newcommand{\Lbar}{\overline{L}}

\renewcommand{\Re}{{\rm{Re}}}

\newenvironment{proof}{
    \noindent{\bf Proof:} \hspace*{1em}}{
    \hspace*{\fill} $\square$\medskip }
\input epsf.tex

\newcommand{\PSbox}[3]{\mbox{\rule{0in}{#3}\hspace{#2}\includegraphics{#1}}}
\begin{document}
\title{Dimers, Tilings and Trees}
\author{ \begin{tabular}{c} Richard W. Kenyon \\
\small CNRS UMR 8628\\  \small Laboratoire de Math\'ematiques\\
\small Universit\'e Paris-Sud\\  \small France\\
\small \texttt{http://topo.math.u-psud.fr/$\sim$kenyon}\\
 \end{tabular} \and \begin{tabular}{c}
\ Scott Sheffield \\
 \small Microsoft Research\\
 \small One Microsoft Way\\
 \small Redmond, Washington\\
 \small U.S.A.\\
 \small\texttt{sheff\char64microsoft.com}
 \end{tabular} }
\date{} \maketitle \begin{abstract} Generalizing results of Temperley \cite{Temp}, Brooks, Smith,
Stone and Tutte \cite{BSST} and others \cite{KPW, traps} we describe a natural equivalence between
three planar objects: weighted bipartite planar graphs; planar Markov chains; and tilings with
convex polygons. This equivalence provides a measure-preserving bijection between dimer coverings
of a weighted bipartite planar graph and spanning trees of the corresponding Markov chain. The
tilings correspond to harmonic functions on the Markov chain and to ``discrete analytic functions''
on the bipartite graph.

The equivalence is extended to infinite periodic graphs, and we
classify the resulting ``almost periodic'' tilings and harmonic
functions.
\end{abstract}

\section{Introduction}
In \cite{Temp}, Temperley gave a bijection between the set of spanning trees of an $n\times n$ grid
and the set of perfect matchings (dimer coverings) of a $(2n-1)\times(2n-1)$ grid with a corner
removed. This bijection was generalized in \cite{KPW} to a weight-preserving bijection (the KPW
construction) from the set of in-directed spanning trees (also known as arborescences) on an
arbitrary weighted, directed planar graph $\Gt$ to the set of perfect matchings on a related graph
$\Gd$. The construction is useful in statistical mechanics because certain types of events in the
spanning tree model can be easily computed using dimer technology, for example winding numbers of
branches and local statistics. For dimer models arising from a spanning tree model, moreover, the
spanning tree formulation provides other useful information. In particular Wilson's algorithm
\cite{Wilson} for generating spanning trees can be used to rapidly simulate dimer configurations.
Moreover, the spanning tree formulation identifies natural boundary conditions (``Temperleyan''
boundary conditions) for the dimer model which allows asymptotic computation of many properties, in
particular conformal invariance properties of dimers \cite{confinv}. However in the paper
\cite{KPW} it was not known if every dimer model on a bipartite planar graph corresponded to a
spanning tree model on a related graph.

A seemingly unrelated construction is the construction of a
``Smith diagram'' from a planar resistor network \cite{BSST}. This
is a tiling of a plane region with squares of arbitrary sizes,
which is associated in a bijective way to a critical-point-free
harmonic function on the network with unit resistances (there is a
square in the tiling for each edge in the graph, whose size is
proportional to the current flow through the edge). This
construction was generalized in \cite{traps} to planar Markov
chains (graphs with transition probabilities), where a harmonic
function gives a tiling with {\it trapezoids}. It was not known at
the time what if any graphical correspondence was natural for {\it
general} polygonal tilings.

In the current paper we extend the above equivalences and describe
a correspondence between these three types of objects: weighted
bipartite planar graphs, planar Markov chains, and tilings with
general convex polygons.

In particular from a weighted bipartite planar graph $\Gd$ we can
construct a tiling $T$ of a plane region with convex polygonal
tiles, and a planar Markov chain $\Gt$, in an essentially
bijective way (that is, up to natural equivalences). There is a
tile in $T$ for each ``white'' vertex of $\Gd$, whose shape is
determined by a {\it discrete analytic function} (see definitions
below). The graph $\Gt$ is a graph on the $1$-skeleton of the
tiles, with transitions determined by their Euclidean geometry.
The tilings are therefore representations of discrete analytic
functions on the bipartite planar graph $\Gd$, which correspond to
harmonic functions on the Markov chain $\Gt$.

An important application of this construction is that it provides
a converse to the Temperley-KPW construction. That is, starting
with the finite weighted bipartite planar graph $\Gd$, one
constructs a Markov chain $\Gt$ and a measure-preserving bijection
from the dimer model on $\Gd$ to the spanning tree process on
$\Gt$. This dimer/spanning tree correspondence has a number of
important consequences. Firstly, it was used in \cite{KOS} in a
fundamental way to classify Gibbs measures on dimer models on
infinite periodic planar graphs. Secondly, since spanning trees
can be sampled efficiently \cite{Wilson}, the construction
provides a way to sample efficiently from bipartite planar dimer
models. Previously the only (provably efficient) way to sample
general planar bipartite dimer models was to do exact computations
of joint edge probabilities.  A third application
\cite{K.surfaces} is that it allows one to compute the asymptotics
of dimer correlations and height fluctuations in terms of the
Green's function on $\Gt$.

In Section \ref{qp} we discuss how the construction extends in the case of infinite periodic
graphs. This is motivated by the study of the dimer model on periodic graphs, see \cite{BP, KOS}.
Given any periodic planar bipartite weighted graph $\Gd$, we produce an essentially unique ``almost
periodic'' planar Markov chain $\Gt$, which extends the dimer/spanning tree correspondence. This
unicity is an important element in the classification theorem of ergodic Gibbs measures on dimer
coverings of $\Gd$ described in \cite{KOS}. \medskip

\noindent{\bf Acknowledgement.}
We thank Andrei Okounkov for several ideas and discussions.

\section{Definitions}
\subsection{Dimers and measures}
Let $\Gd=(V,E)$ be a finite bipartite planar graph. Bipartite
means that the vertices $V$ can be $2$-colored, that is, colored
black and white so that black vertices are only adjacent to white
vertices and vice versa.   Let $\nu\colon E\to(0,\infty)$ be a
weight function on the edges. A {\bf perfect matching}, or {\bf
dimer configuration} $M\subset E$ is a set of edges with the
property that each vertex is contained in exactly one edge in $M$.
The weight of a matching $M$ is $\nu(M)=\prod_{e\in M} \nu(e)$.
Let $\M(\Gd)$ denote the set of perfect matchings of $\Gd$. Let
$\mu$ be the probability measure on $\M(\Gd)$ giving a matching a
probability proportional to its weight: $\mu(M)=\frac1{Z}\nu(M)$
where $Z=\sum_{M\in\M(\Gd)}\nu(M).$

\subsection{Kasteleyn matrices}
If $\Gd$ has $n$ black and $n$ white vertices, a {\bf Kasteleyn
matrix} (see \cite{Kast}) for $\Gd$ is a real $n \times n$ matrix
$K=(K_{i,j})$ whose rows index the black vertices and columns
index the white vertices of $\Gd$, defined as follows. The entry
$K_{i,j}$ is zero if there is no edge from $b_i$ to $w_j$, and if

there is an edge of weight $\nu(b_iw_j)$ then $K_{i,j} =\pm
\nu(b_iw_j)$, where the signs are chosen so that the product of
signs of edges around every interior face of $K$ is
$(-1)^{d/2+1}$, where $d$ is the degree of the face.  This
property of signs is not changed if we multiply all elements in a
particular column or row of $K$ by $-1$ (because each vertex of
$\Gd$ has an even number---zero or two---of edges on each face of
$\Gd$).  Moreover, such a choice of signs always exists, and by
Kasteleyn's theorem, the determinant of $K$ is (up to sign) the
sum of the weights of the matchings of $\Gd$ (\cite{Kast}). By a
{\bf discrete analytic function} we mean a function $f$ on black
vertices (resp. white vertices) which satisfies $fK=0$ (resp.
$Kf=0$). This generalizes the definition of discrete analytic
function on $\Z^2$ defined in \cite{Duffin, confinv}. These
functions play a role implicitly in sections
\ref{constructionsection} and \ref{qp}.

\subsection{Gauge transformations}
If we multiply the weights of all the edges in $\Gd$ having a
fixed vertex by a constant, the measure $\mu$ does not change,
since exactly one of these weights is used in every configuration.
More generally, two weight functions $\nu_1,\nu_2$ are said to be
{\bf gauge equivalent} if $\nu_1/\nu_2$ is a product of such
operations, that is, if there are functions $f_1$ on white
vertices and $f_2$ on black vertices so that for each edge $wb$,
$\nu_1(wb)/\nu_2(wb)=f_1(w)f_2(b).$ Gauge equivalent weights
define the same measure $\mu$.

Multiplying the $i$th row (resp., column) of a Kasteleyn matrix
$K$ by a positive, non-zero constant $c$ is equivalent to
multiplying by $c$ the weights of all of the edges of $\Gd$
incident to $b_i$ (resp., $w_i$). In other words any matrix
$\tilde K$ obtained from $K$ by multiplying the rows and columns
of $K$ by non-zero constants will be a Kastelyn matrix for a graph
which is gauge equivalent to $\Gd$.

\section{T-graphs and corresponding dimer/spanning-forest models}
In this section, we define a family of planar graphs called
T-graphs and describe a weight-preserving correspondence between
the spanning trees on a T-graph and dimer configurations on a
bipartite graph derived from it. This is closely related to the
result of \cite{KPW}, who also give a relation between spanning
trees on a general planar graph and dimers on a derived graph.

In the present context however our derivation can be reversed: we
will see in Section \ref{constructionsection} that for every
bipartite planar graph, which is {\it non-degenerate} in the sense
that it contains no edges which fail to be used in any perfect
matching of the graph (for the purposes of the dimer model, it
makes sense to delete these edges), endowed with a generic choice
of weights, there is a gauge-equivalent graph which can be derived
from a T-graph in this way.  By taking limits, the correspondence
generalizes to the case when the weights are not assumed to be
generic.

\subsection{Complete edges that form T-graphs}
The definition of $T$-graph on a torus---which we use in section \ref{qp}---is quite simple. A
disjoint collection $L = \{ L_1, L_2, \ldots, L_n \}$ of open line segments in the torus $\R^2 /
\Z^2$ {\bf forms a T-graph in the torus} if $\cup_{i=1}^n L_i$ is closed. The term ``T-graph''
refers to the fact each endpoint of a given $L_i$ necessarily lies on the interior of some $L_j$
with $j \not = i$. In other words, each $L_i$ ``tees into'' an $L_j$ at each of its two endpoints.

\begin{figure}\begin{center}\leavevmode \epsfbox[57 608 415 730]{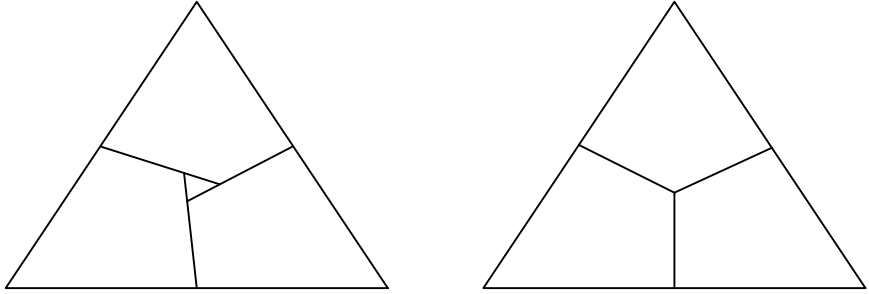}\end{center}
\caption{\label{nongeneric}(a) Line segments that form a
$T$-graph. (b) Line segments that don't form a T-graph.}
\end{figure}

We say a disjoint collection $L_1, L_2, \ldots, L_n$ of open line
segments in $\R^2$ {\bf forms a T-graph in $\R^2$} if
$\cup_{i=1}^n L_i$ is connected and contains all of its limit
points except for some set $R = \{r_1, \ldots, r_m \}$, where each
$r_i$ lies on the boundary of the infinite component of $\R^2$
minus the closure $\Lbar$ of $\cup_{i=1}^n L_i$.  Elements in $R$
are called {\bf root vertices}. For example, a single open line
segment forms a T-graph with root vertices given by the two
endpoints.  A pair of open line segments---one of which tees into
the other to make the letter ``T''---forms a T graph with three
root vertices. The three open edges of a triangle also form a
T-graph with three root vertices. A partitioning of a convex
polygon $P$ into convex polygonal tiles using a finite number of
line segments will form a $T$-graph with root vertices at the
vertices of $P$ if and only if it is {\bf generic} in the sense
that the endpoint of each of these line segments lies either on
the interior of another line segment or on the boundary of $P$.
(See Figure \ref{nongeneric}.)

Note that each endpoint of a given $L_i$ is {\bf either} a root vertex or an interior point of some
$L_j$.  To distinguish the $L_i$ from subsegments of the $L_i$ (which we discuss later) we refer to
the $L_i$ as {\bf complete edges}.

In subsequent subsections, we will use $L$ to define a weighted, directed graph $\Gt(L)$ and a
weighted, bipartite graph $\Gd(L)$.  The ultimate goal of this section will be to derive a
weight-preserving bijection between directed spanning forests on $\Gt(L)$ (with specified roots)
and perfect matchings of $\Gd(L)$.  When the choice of $L$ is clear from the context, we write $\Gt
= \Gt(L)$ and $\Gd = \Gd(L)$.

\subsection{T-graphs and their duals} \label{tgraphdualsubsection} The set $V_T(L)$ of vertices of $\Gt$
(the graph we will call the {\bf tree-graph} of $L$) is the set of points in $\R^2$ which are
endpoints of at least one of the $L_i$.  A vertex $v$ which is in the interior of a complete edge
$L_i$ (called an {\bf interior vertex}) has exactly two edges in $\Gt$ directed outwards from it:
these edges point towards the two immediate neighbors, $v_1$ and $v_2$, along $L_i$ (one on each
side of $v$). The weights on the edge from $v$ to these two $v_i$ are chosen in such a way that the
two weights add up to one and are inversely proportional to the Euclidean distances $|v - v_i|$.
These weights correspond to the transition probabilities of a Markov chain on $V_T(L)$. The root
vertices are sinks of $\Gt$ (they have no outgoing edges in $\Gt$) and are fixed points of the
Markov chain. Note that (by our choice of transition probabilities) the expected change in
Euclidean position during a step of the Markov chain is always zero; thus, a random walk on $\Gt$
--- viewed as a Markov chain on positions in $\R^2$ --- is a martingale.  In other words, the
coordinate functions on the vertices of $\Gt$ are harmonic functions on $\Gt$ away from the root
vertices.

\begin{figure}[htbp] \PSbox{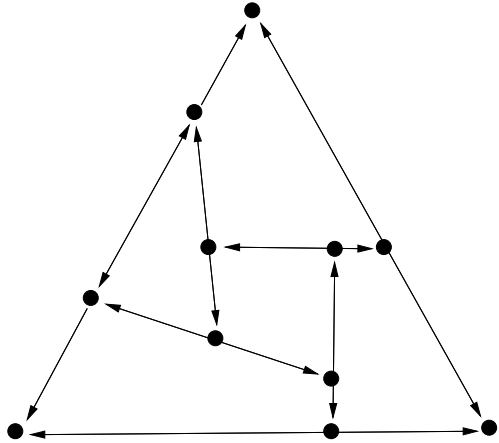}{0in}{2.5in}
\caption{The directed Tree-graph $\Gt$. The
root vertices are the corners of the triangle. \label{Gb}}
\end{figure}

See Figure \ref{Gb} for an example of a T-graph with three roots.  Note that, by convention, when
we have transitions both from $i$ to $j$ and from $j$ to $i$, rather than drawing two directed
edges in the graph $\Gt$ we draw a single edge with two transition probabilities, one from each
end.

We define $\Gt' = \Gt'(L)$, an undirected dual graph of $\Gt$, as follows. Let $C$ be an arbitrary
simple closed curve that encircles the $\cup_{i=1}^n L_i$ and contains each of the root vertices
$r_1, \ldots, r_m$ in clockwise order. The vertices of $\Gt'$ are the bounded faces of $\Gt\cup C$
(bounded connected components of $\R^2\backslash(\cup_{i=1}^n L_i\cup C)$). Faces of $\Gt\cup C$
adjacent to $C$ are called {\bf outer faces} of $\Gt$: they correspond to {\bf outer vertices} of
$\Gt'$. Two vertices of $\Gt'$ are connected by an edge $a$ of $\Gt'$ if the corresponding faces of
$\Gt$ are adjacent across an edge of $\Gt$. For an edge $e$ of $\Gt$ we denote by $e^*$ its
corresponding dual edge.

\begin{lem} \label{vertnum} If $L_1,L_2,\ldots, L_n$ form a T-graph, then $\Gt$ has exactly $n+1$
faces (including outer faces).  Hence $\Gt'$ has $n+1$ vertices. \end{lem}

\begin{proof} This follows from Euler's formula.  The line segments $L_i$ decompose the interior of
$C$ into some number $n_2$ of open faces (open $2$-cells), $n_1 = n$ open complete edges, and $n_0
= 0$ vertices. Since $n_2-n_1+n_0=1$, the Euler characteristic of the disc, the result follows.
\end{proof}

\subsection{Spanning trees}
A {\bf spanning tree} of a graph $G$ is a subset of edges which is connected, contains no cycle,
and passes through every vertex. If the edges of $G$ are directed, a {\bf directed spanning tree},
or arborescence, is a spanning tree in which every vertex but one (called the root vertex) has a
unique outgoing edge. Given a subset of vertices of $G$ called root vertices, a {\bf directed
spanning forest} is a set of edges with no cycles, passing through all vertices, each non-root
vertex having a unique outgoing edge, and each component of which is connected to a unique root
vertex.

\begin{figure} \begin{center} \leavevmode \epsfbox[50 421 365 743]{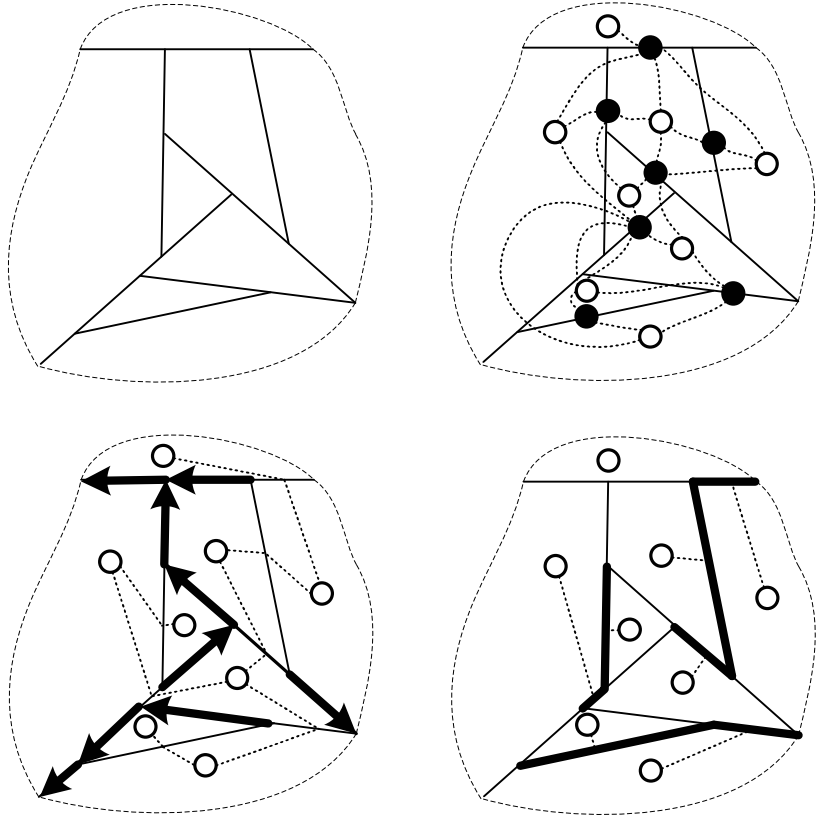} \end{center}
\caption {(a) Edges of a T-graph $L$ and the surrounding curve $C$. (b) The graph $\Gdr$.  This
graph is the incidence graph for the set of complete edges of $L$ (drawn as black vertices) and
faces of $L$ (drawn as white vertices). (c) A spanning forest of $\F$ of $\Gt$, drawn with thick
arrows, and the dual spanning tree $\Tr$ on $\Gt'$, drawn with dotted lines connecting vertices of
$\Gt'$ (which are faces of $\Gt$, represented as white vertices).  Each such dotted line crosses a
segment of a complete edge that is {\it not} used in $\F$ (there is exactly one such segment for
each complete edge).  (d) The marked matching corresponding to $\F$ when the dual root is taken to
be the (unmatched) uppermost white vertex.}\label{tiling} \end{figure}

We will employ the following correspondence between (non-directed) spanning trees in $\Gt'$ and
their (non-directed) {\bf dual spanning forests} in $\Gt$.  Using the correspondence between edges
of $\Gt'$ and edges of $\Gt$, we can think of edge subsets of both $\Gt$ and $\Gt'$ as subsets of
the set of all edges of $\Gt$. Using this interpretation, we state the following lemma (which is
illustrated in Figure \ref{tiling}):

\begin{lem} \label{treetodual} The complement of a spanning tree $\Tr$ of $\Gt'$ is a spanning
forest $\F$ of $\Gt$, with roots at the root vertices.  Similarly, the complement of a spanning
forest $\F$ of $\Gt$, with roots at the root vertices, is a spanning tree $\Tr$ of $\Gt'$.
\end{lem} \begin{proof} We sketch the standard tree dualization argument. If $\F$ is a spanning
forest of $\Gt$, with roots at root vertices, its complement $\Tr$ cannot contain any cycles in
$\Gt'$ (since such a cycle would separate at least one interior vertex of $\Gt$ from the root
vertices), and it must be connected (since otherwise, the set of edges separating two components
would either form a cycle in $\Gt$ or a path connecting two root vertices in $\Gt$); hence it is a
spanning tree. Similarly, if $\Tr$ is a spanning tree of $\Gt'$, its complement $\F$ cannot contain
cycles of $\Gt$ (since such a cycle would separate at least one inner face of $\Gt$ from the outer
faces) and each connected component of $\F$ contains at least one root vertex (since otherwise the
set of edges separating that component of $\F$ from its complement would form a cycle in $\Gt$).
\end{proof}

\subsection{Dimer graphs from $T$-graphs} \label{dimercorrespondencesubsection}
Now we will define the weighted, bipartite (non-directed) graph $\Gd=\Gd(L)$.  First, we define a
slightly larger $\Gdr = \Gdr(L)$, whose black vertices are the $n$ complete edges $L_i$ and whose
white vertices are the $n+1$ faces of $\Gt$ (including outer faces).

A white vertex $w$ of $\Gdr$ is adjacent to a black vertex $b$ of $\Gdr$ if the face $F$
corresponding to $w$ contains a portion of the $L_i$ corresponding to $b$ as its boundary.  The
weight $\nu((w,b))$ is then given by the Euclidean length of the portion of the line segment.  The
graph thus defined is planar. To see this, note that it can be drawn on top of the tiling $\Lbar$
as follows: put a white vertex in the interior of each face, and a black vertex in the center of
each complete edge. When $w$ and $b$ are connected, draw a line from $w$ inside the corresponding
face towards the complete edge corresponding to $b$, and then along this complete edge, staying
just to one side, until the center is reached.  It is not hard to see that this can be done in such
a way that the paths do not intersect.

The graph $\Gd$ is formed from $\Gdr$ by (arbitrarily) picking one of the outer white vertices of
$\Gdr$ and removing it; we will refer to the removed vertex as the {\bf dual root} of $\Gd$. Now,
$\Gd$ is a weighted bipartite graph with $n$ white and $n$ black vertices. To every edge $e =
(w,b)$ in a perfect matching of $\Gd$ (where $w$ corresponds to a face $F$ and $b$ to a complete
edge $L_i$), we denote by $S_e$ the segment of the $L_i$ which borders $F$.  Because of our choice
of weights, $\mu(M)$ (where $\mu$ is the probability measure on perfect matchings defined in the
introduction) is proportional to $\prod_{e \in M} |S_e|$ where $|S_e|$ is the Euclidean length of
$S_e$.  Now, the edge segment $S_e$ may have vertices of $\Gt$ in its interior; these vertices
divide $S_e$ into subsegments, each of which has vertices of $\Gt$ as its endpoints and hence
corresponds to an edge of $\Gt$.  A {\bf marked matching} of $\Gd$ is a matching $M$ of $\Gd$
together with a specified subsegment $S_e'$ of $S_e$ (which, again, we may interpret as an edge in
$\Gt$) for each $e \in M$. We extend $\mu$ to give a measure on random marked matchings as follows:
to sample a random marked matching, first choose a random matching. Then for each edge $e$, choose
an $S_e'$ from among the subsegments of $S_e$, where probability of each subsegment is proportional
to its length. If $M'$ is a marked matching, then $\mu(M')$ is proportional to $\prod_{e \in M}
|S_e'|$.

\begin{center} \begin{table}[ht] \label{graphtable}
\begin{tabular}{ l l } \textbf{Graph}    &  \textbf{Vertex Set} \\
\hline
$\Gt = \Gt(L)$            &  Points that are endpoints of some $L_i$.  \\
$\Gt'$       &  Faces and outer faces of
$\Gt$ (i.e., bounded components of $\R^2 \backslash (L \cup C)$). \\
$\Gdr$    &  Faces and outer faces of $\Gt$ (one partite class).\\
    &        Complete edges in $L$ (other partite class).\\
$\Gd$     &  Same as $\Gdr$ but with one outer face (the dual root) omitted. \\
$\Gd'$     &  Faces and some outer faces of $\Gd$ (which correspond to vertices of $\Gt$).  \\
\hline \end{tabular} \caption{Summary of graphs constructed from a collection of edges $L$ that
forms a $T$-graph.  The graph $\Gd'$ (which is not exactly the same as $\Gt$) is defined precisely
in Section \ref{Kastelynflowsubsection}.} \end{table} \end{center}

\subsection{From dimers to trees}
Let $\Tr_{M'} = \{ S_e': e \in M \}$ be the set of of edges corresponding to a marked matching $M'$
of $\Gd$.  Each $S_e'$ corresponds to an edge of $\Gt'$, so we can think of $\Tr_{M'}$ as a
subgraph of $\Gt'$.  We direct each such edge of $\Tr_{M'}$ (corresponding to some $S_e'$) of this
graph from the face which corresponds to a vertex in $e$ towards the face which does not. We use
this interpretation of $\Tr_{M'}$ (as a directed subgraph of $\Gt'$) in the following lemma:

\begin{lem} \label{matchingtotree} If $M'$ is a marked matching, then $\Tr_{M'}$ is an in-directed
spanning tree of $\Gt'$, rooted at the dual root.  The dual $\F_{M'}$ of $\Tr_{M'}$ is thus a
spanning forest of $\Gt$ (when $\Gt$ is viewed as an undirected graph). \end{lem}

\begin{proof} It is sufficient to prove that $\Tr_{M'}$ has no directed cycles, since it contains
exactly one edge pointing away from each face of $\Gt$ (excluding the dual root). This is
accomplished using Euler's formula. Suppose that $\Tr_{M'}$ had a directed cycle $F_0, F_1, \ldots,
F_j = F_0$ of faces of $\Gt$. Let $S_i$ be the segment $S_e'$ separating $F_i$ and $F_{i+1}$. Let
$C'$ be a simple closed curve which starts in the interior of $F_0$, passes through $S_0$ at one
point, moves through the interior of $F_1$, passes through $S_1$ at a single point, etc. until it
returns to $F_0$. Except for its intersections with the $S_i$'s, each at a single point, $C'$ is
entirely contained in the union of the interiors of the $F_i$. The intersection of the $L_i$ with
the interior of $C'$ gives a decomposition of this interior into $n_2$ two-cells (where $n_2$ is
the number of faces partial or completely contained inside the loop $C'$), $n_1$ open one-cells and
$n_0=0$ vertices. Thus, by Euler's formula $n_2 - n_1 + n_0 = n_2-n_1=1$. In particular $n_2+n_1$
is odd.

However, the sequence $w_0, b_0, w_1, b_1, \ldots, b_{j-1}, w_j = w_0$ (where $w_i$ is the white
vertex of $\Gd$ corresponding to $F_i$ and $b_i$ is the black vertex corresponding to the complete
edge containing $S_i$) is a cycle in $\Gd$, alternating edges of which are contained in $M$. The
set of vertices in $\Gd$ enclosed by this cycle must be matched only with each other in a perfect
matching (since the cycle disconnects these vertices from the rest of the graph). This is a
contradiction to the fact that $n_2+n_1$ is odd. \end{proof}

Let $\mu_F$ be the measure on directed spanning forests of $\Gt$, rooted at the root vertices, for
which $\mu_F(\F)$ is proportional to the product of the weights of the edges in $\F$.  Since each
of the two outgoing edges of a given interior vertex has weight (by construction) inversely
proportional to its Euclidean length, $\mu_F(\F)$ is inversely proportional to the product of the
lengths of the edges of $\F$; hence, $\mu_F(\F)$ is also proportional to the product of the
Euclidean lengths of all edges of $\Gt$ which do not appear (directed one way or another) in $\F$.

The following is the main result of this section.

\begin{thm} \label{corrthm} The map $M' \rightarrow \F_{M'}$ gives a one-to-one correspondence
between marked matchings of $\Gd$ and in-directed spanning forests of $\Gt$, rooted at $R$. The
correspondence is measure preserving, i.e., $\mu(M') = \mu_F(\F_{M'})$. \end{thm}

\begin{proof}  First, we would like to interpret $\F_{M'}$ as a directed spanning forest of $\Gt$
by orienting each edge of $\F_{M'}$ towards its root vertex.   In order to do this, we must check
that if $M'$ is a perfect matching of $\Gd$, then the directed path along $\F_{M'}$, from a vertex
$v$ to a root vertex, is a directed path of $\Gt$.  To see this, note first $\F_{M'}$ contains all
but one segment of each of the $L_i$; thus, for every interior vertex $v$ of $\Gt$ (interior to
some $L_i$), $\F_{M'}$ includes a path from $v$ to exactly one of the endpoints of $L_i$.  Call
this vertex $v_1$; each of the directed edges in the directed path from $v$ to $v_1$ is a directed
edge of $\Gt$.  If $v_1$ is also an interior vertex of some $L_j$, then there is a path of edges in
$\F_{M'}$ from $v_1$ to some endpoint $v_2$ of $L_j$.  Iterating this process, we must eventually
produce a directed path from $v$ to a root (since $\F_{M'}$ has no cycles).

It now follows immediately from Lemma \ref{matchingtotree} and our choice of weights, that
$\mu_F(\F_{M'})$ is proportional to $\mu(M')$, since each is proportional to the same product of
edge lengths.  The proof we will be complete once we show that the map $M' \rightarrow \F_{M'}$ is
invertible.

Let $\F$ be an arbitrary directed spanning forest $\F$, rooted in $R$.  Since only the endpoints of
a given $L_i$ have outgoing edges pointing to vertices not on $L_i$, each vertex of $L_i$ belongs
to a path pointing to one of the two endpoints.  It follows that $\F$ must include all but one of
the subsegments of $L_i$.  By Lemma \ref{treetodual}, the dual of $\F$ is a spanning tree of
$\Gt'$, which we may view as being directed towards the dual root.  Each face $F$ of $\Gt$ is
(besides the dual root) is directed towards another face across an edge segment of one of the
$L_i$.  Pairing of $F$ with the edge segment produced in this way gives a marked matching $M'$ for
which $\F_{M'} = \F$. \end{proof}

\subsection{T-graphs and dimers on the torus}
If $L = \{L_1, \ldots, L_n\}$ forms a T-graph on the torus, then we can construct $\Gt = \Gt(L)$
exactly as above; in this case, $\Gt(L)$ has no root vertices and no outer faces.  Since the faces
of $\Gt$ and open edges $L_i$ give a decomposition of the torus into one-cells and two-cells,
Euler's formula implies that $\Gt$ has exactly $n$ faces.  We construct $\Gd$ as above (with white
vertices given by faces $F$ of $\Gt(L)$, black vertices by the complete edges $L_i$, and edges
occurring between $F$ and $L_i$ that share a line segment, weighted according to the length of that
segment). We also construct $\Gt'$ in a similar fashion.

A {\bf cycle-rooted spanning forest} $\F$ of $\Gt$ is a (directed) subgraph of $\Gt$---with one
outgoing edge from each vertex of $\Gt$---which has no null-homotopic (directed) cycles (i.e., no
cycles which---when lifted to the universal cover of the torus---start and end at the same place).
The ``roots'' of such an $\F$ are the directed cycles of $\F$.  Clearly, every such $\F$ has at
least one (non-null-homotopic) directed cycle.

The dual of $\F$ is a cycle-rooted spanning forest $\F'$ on $\Gt'$.  Now, if $\F$ has exactly $j$
cycles, then it is not hard to see that $\F'$ has $j$ cycles as well.  We can view $\F'$ as a
directed cycle-rooted spanning forest by directing each edge not on a cycle towards its cycle root;
and then orienting all of the edges in a given cycle one of the two possible directions (there are
$2^j$ ways of doing this).  The proof of the following is now similar to the proof of Theorem
\ref{corrthm}.

\begin{thm} There is a one-to-one weight preserving correspondence between perfect matchings on
$\Gd$ and in-directed cycle-rooted spanning forests $\F'$ on $\Gt'$ whose dual cycle-rooted
spanning forests $\F$ are in-directed, cycle-rooted spanning forests of $\Gt$. \end{thm}

\begin{figure}[htbp]
\vskip3in
\begin{center} \leavevmode \epsfbox[250 550 300
600]{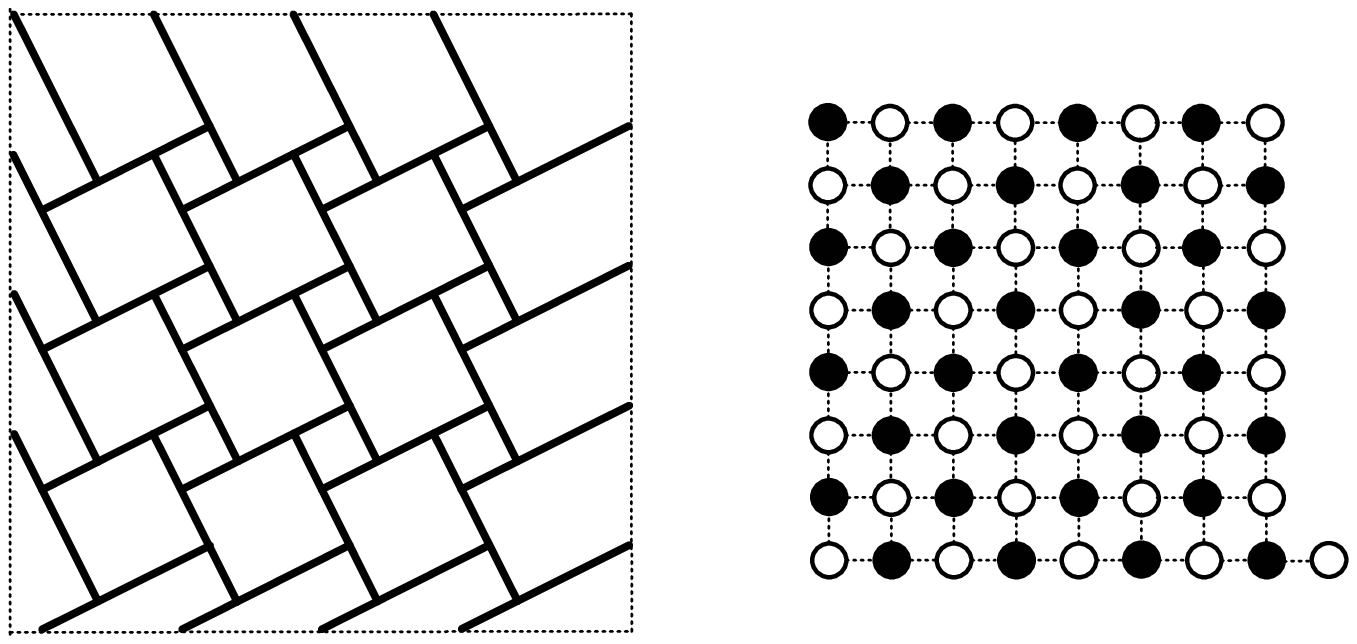}
\end{center} \caption {A T-graph $\Gt$ in the plane and the corresponding graph
$\Gdr$.}\label{periodicexample} \end{figure}

T-graphs in a torus can be extended to give periodic T-graphs on the plane, finite subsets of which
correspond to finite subgraphs of infinite lattice graphs, such as the grid graph in Example
\ref{periodicexample}.

\section{T-graphs from dimer graphs} \label{constructionsection} In this section, we describe a procedure
for generating $\Gt$ from $\Gd$ that applies whenever the so-called {\bf Kastelyn matrix} fails to
have certain degeneracies.  Before we begin the construction, we will define Kastelyn matrices and
say a word about the kinds of graphs for which these degeneracies occur.

\subsection{Cuts, and breakers}

Say a square matrix $K$ is {\bf $k$-degenerate} if it has an
$(n-k) \times (n-k)$ minor whose determinant is zero; otherwise it
is {\bf $k$-non-degenerate}.  The following lemma follows from the
standard correspondence between determinants of $k$ minors of
$K^{-1}$ and $(n-k)$ minors of $K$:
\begin{lem} $K$ is
$0$-non-degenerate if and only if it is invertible. Assuming $K$
is invertible, $K$ is $k$-non-degenerate if and only if $K^{-1}$
is $(n-k)$-non-degenerate. \end{lem}

Suppose now $K$ is a Kasteleyn matrix for a bipartite planar graph
$\Gd$. The following is immediate:

\begin{lem} If $K$ and $\tilde K$ are gauge equivalent, then $K$ is $k$-degenerate if and only if
$\tilde K$ is $k$ degenerate \end{lem}

A bipartite graph is {\bf balanced} if it contains an equal number of black and white vertices. A
{\bf $k$-cut} $A$ of a balanced bipartite graph $\Gd$ is subset of the vertices for which:
\begin{enumerate}
\item $A$ contains at least one white vertex
\item $A$ contains $k$ more black vertices than white vertices \item Each edge of $\Gd$ that
connects $A$ to its
complement has a black vertex in $A$. \end{enumerate} Note that if $A$ is a $k$-cut, then its
complement would be a $k$-cut if the colors black and white were reversed.  In particular, the
existence of $k$ cuts does not depend on which of the two ways we choose to color the vertices.
Also, if $A$ is a $k$-cut of $\Gd$, then by adding black vertices to $A$ and/or removing white
vertices from $A$, we can construct $m$-cuts for any $k \leq m \leq n-1$.  An obvious parity
argument implies the following:

\begin{lem} If $A$ is a $k$-cut of $\Gd$, then any perfect matching of $\Gd$ contains exactly $k$
edges which connect $A$ to its complement; each of these edges matches a black vertex of $A$ and a
white vertex of its complement. \end{lem}

A {\bf $k$-breaker} is a subset $S$ of the vertices of $\Gd$ with exactly $k$ white and $k$ black
vertices for which the induced subgraph $\Gd \backslash S$ of $\Gd$ has no perfect matchings.

\begin{lem} If $\Gd$ is a connected, balanced, bipartite graph, then $\Gd$ \begin{enumerate} \item
has a $(-1)$-cut if and only if it has no perfect matching. \item has a $0$-cut if and only if $\Gd$
contains {\bf unused edges} (i.e., edges which occur in no perfect matching of $\Gd$). \item
generally has a $k$-cut if and only if has a $(k+1)$-breaker. \end{enumerate} \end{lem}

\begin{proof} The first item is an immediate consequence of the Hall marriage theorem.  That
theorem states that $\Gd$ has a perfect matching if and only if there is no set $B$ such that $B$
has $m$ more white vertices than black vertices and there are fewer than $m$ edges connecting white
vertices of $B$ to its complement.   A $(-1)$-cut is clearly such a set, with $m=1$. Conversely,
given $B$ as described above, construct $B'$ by removing from $B$ all of the (at most $m-1$) white
vertices of $B$ connected to the complement of $B$, and if necessary, some arbitrary additional
white vertices (so that $m-1$ vertices are removed in all).  Then $B'$ is a $(-1)$-cut.

For the second item, first, it is clear that if $A$ is a $0$-cut of $B$, then all of the edges
connecting $A$ to its complement will be unused. Conversely, if $\Gd$ has an unused edge $e$, then
the graph $G$ formed by removing edge $e$ and its two vertices from $\Gd$ will not have any perfect
matching. Therefore it will have a $(-1)$-cut $A$ by part 1.  The union of $A$ and the black vertex
of $e$ is a thus a $0$-cut.  (Aside: if $\Gd$ has a {\bf forced edge}---i.e., an edge $e$ which
occurs in {\it every} perfect matching of $\Gd$---then all the edges that share vertices with $e$
will be unused.)

The same argument implies the third statement in the case $k=0$. For larger $k$, if $\Gd$ has a
$k$-cut $A$, then any subset of $(k+1)$ black vertices of $A$ and $(k+1)$ white vertices of its
complement is a $(k+1)$-breaker (since the remaining set of vertices in $A$ contains more white
than black vertices, but there are no edges connecting white vertices of this remaining set with
its complement). Conversely, if $S$ is a $(k+1)$-breaker, then $\Gd \backslash S$ has a $(-1)$-cut
$A$, and the union $A$ and the black vertices of $S$ is a $k$-cut of $\Gd$. \end{proof}

\begin{lem} \label{degeneratelemma} $\Gd$ has no $k$-breaker (or, equivalently, no $(k-1)$-cut) if
and only if, for a generic choice of positive weights of the edges of $\Gd$, the Kastelyn matrix $K
= K(\Gd)$ is $k$-non-degenerate. \end{lem} \begin{proof} The determinant of an $(n-k) \times (n-k)$
minor of the Kastelyn matrix is a polynomial of the edge weights. Clearly, this polynomial will be
zero for a given minor precisely when the set of $k$ white and $k$ black vertices corresponding to
rows and columns not in the minor is a $k$-breaker.  The result follows from the fact that any
non-zero polynomial in finitely many real variables is non-zero for a generic choice of inputs.
\end{proof}

In this paper, we will mainly be interested in whether $K$ is $k$-degenerate for $k \in \{0,1,2\}$.
But we know that whenever $K$ is a Kastelyn matrix (of a graph having a perfect matching) it is
$k$-non-degenerate for $k=0$.  And assuming $\Gd$ has no unused edges (which we may always assume
throughout, since the perfect matching model will be unchanged if we remove unused edges from
$\Gd$) $K$ is generically $1$-non-degenerate.  We will address the potential failure of $K$ to be
$2$-non-degenerate in a later section.

\subsection{T-graphs: construction via integration of Kastelyn flow} \label{Kastelynflowsubsection}
Let $\Gd$ be a finite, weighted bipartite planar graph (with
positive generic weight function $\nu$) with $n$ black vertices $b_1, b_2, \ldots, b_n$ and $n$
white vertices $w_1, w_2, \ldots, w_n$.  Suppose $\Gd$ has a perfect matching and no unused edges.
Suppose that $\Gd$ has no $1$-cuts---and hence each of the entries and two-by-two minors of
$K^{-1}$ is non-zero (i.e., $K$ is $1$-non-degenerate and $2$-non-degenerate).

We will now construct a T-graph corresponding to $\Gd$ in the case that $K$ is $2$-non-degenerate.

First, we may think of $K$ as describing a linear map from the space $\R^W$ of functions on white
vertices to the space $\R^B$ of functions on black vertices. Let $b_0$ be a fixed vertex on the
outer boundary of $\Gd$. Suppose that $\Gd$ has $m$ black and $m$ white vertices on its outer
boundary face.  Fix a generic convex $m+1$-gon $Q$ with edge vectors $q_0,\dots,q_m\in\C$ in cyclic
order (and $q_0=-\sum_{i=1}^m q_i$). Vertices of $Q$ will be the root vertices of $\Gt$.  Suppose
that $A_w \in \R^W$ assumes the values $q_1, \dots, q_m$ in cyclic order on the white vertices on
the boundary face, and that $A_w$ vanishes on all other white vertices of $\Gd$. Let $A_b$ be the
function on black vertices which is equal to $1$ at $b_0$ and $0$ everywhere else. Denote by $\1$
the all-ones column vector and by $\1^t$ its transpose.  View $A_b$ as a column vector and $A_w$ as
a row vector.

We claim that there is a unique matrix $\tilde K$, gauge equivalent to $K$, for which $\tilde K \1$
is a non-zero multiple of $A_b$ and $\1^t \tilde K = A_w$. The matrix $\tilde K$ can be derived
explicitly from $K$ as follows. Since $K$ is invertible, there exists a vector $f$ for which $K f =
A_b$. Multiplying the $i$th column of $K$ by the $i$th component of $f$ (non-zero, because $K$ is
$1$-non-degenerate) produces a $K'$ for which $K' \1 = A_b$.  Next, there exists a row vector $g$
for which $g K' = A_w$. Multiplying the $j$th row of $K'$ by the $j$th component of $g$ (also
non-zero, since $(K')^{-1}$ is $1$-non-degenerate and nonzero entries of $A_w$ are generic) gives
the desired $\tilde K$.

We may think of $\tilde K$ as describing a vector flow
($2$-component flow) on $\Gd$: sending $\tilde K_{i,j}$ units of
flow from $b_i$ to $w_j$.  The net flow into each non-boundary
white vertex and each black vertex (except $b_0$) is zero.  Now,
draw a dotted line from each white vertex on the outer face of
$\Gd$ to infinity, and from $b_0$ to infinity, so as to divide the
outer face of $\Gd$ into $m+1$ outer faces; take these faces and
the interior faces of $\Gd$ as the vertices of the dual graph
$\Gd'$ of $\Gd$.  Then $\tilde K$ also describes a dual flow on
$\Gd'$ (obtained by rotating each edge ninety degrees
counter-clockwise) whose net flow around each non-boundary face of
$\Gd'$ is zero; viewed in this light, $\tilde K$ is the gradient
of a function $\psi: \Gd' \rightarrow \C$.

Now, we claim that each pair of (complex) components of $g$ is linearly independent (as a pair of
vectors in $\C =\R^2$). To see this, let $a$ and $b$ be basis column vectors, so that
$(K')^{-1}(a)$ and $(K')^{-1}b$ are columns of the matrix $(K')^{-1}$. Since the determinants of
the two-by-two minors of $(K')^{-1}$ are non-zero, no complex component of the vector $z =
(K')^{-1}a + i(K')^{-1}b = (K')^{-1}(a + ib)$ is a real multiple of any other component of that
vector (in particular, all of the components of $z$ are non-zero). Now $A_w$ is a generic linear
combination of vectors of the above form $a+ib$, so no component of $g = K^{-1}(A_w)$ is a real
multiple of any other component of $g$.

Since $K'$ is real, all the components of $\tilde K$ in a given row are nonzero complex numbers
lying on the same line through the origin, and the  directions are different in each row.

Now, extend $\psi$ linearly to the edges of $\Gd'$, so that $\psi$ maps each edge to a line
segment.  For each black vertex $b_i$ of $\Gd$, corresponding to a black face of $\Gd'$, the $\psi$
image of the union of the edges incident to $b_i$ is a line segment, whose interior we denote by
$L_i$; the above argument implies that no two of the $L_i$ are parallel.

Here is the main result. \begin{thm} \label{TgraphfromK} If $K$ is $2$-non-degenerate then the $L =
\{ L_1, \ldots, L_n \}$ defined above forms a T-graph with root vertices at vertices of $Q$, and
$\Gd = \Gd(L)$ (up to gauge equivalence). Moreover, if $v$ is a vertex of $\Gd'$, then $\psi(v)$ is
a vertex of the T-graph; if $v$ corresponds to an outer face of $\Gd$, then $\psi(v)$ is a root
vertex of the T-graph. \end{thm}

\begin{proof} First, the change in $\psi$, as one moves from outer face $F$ of $\Gd$ around a
vertex $v$ to another outer face, is given by the flow of $\tilde K$ into $v$, which is given by
$q_i$, the $i$th component of $A_w$, whenever $v$ is a white vertex $w_i$, and zero when $v$ is any
black vertex besides $b_0$. By moving around the polygon in steps, it is clear that (up to an
additive constant) $\psi(F)$ assumes the values of the vertices of the convex polygon in cyclic
order.

Let $f$ be an interior vertex of $\Gd'$.  We claim that for some black face incident to $f$, with
vertices $f_1$ and $f_2$ incident to $f$, $\psi(f_1) - \psi(f)$ and $\psi(f_2) - \psi(f)$ point in
opposite directions.  Suppose otherwise.  Then $\tilde K$ would have to assume opposite signs on
the entry corresponding to each such pair of edges $(f,f_1)$ and $(f,f_2)$.  By the definition of a
Kastelyn matrix, $\tilde K$ has positive sign for an odd (resp., even) number of the edges incident
to $f$ if the total number of edges is $0\bmod 4$ (resp., $2\bmod 4$), so this is a contradiction.
It follows that $\psi(f)$ is an interior vertex of at least one $L_i$. In particular, this implies
that the endpoint of each $L_i$ is either an interior vertex of some $L_j$ or a root vertex.

It also implies a {\bf maximal principle}, i.e., that for any vector $u$ in $\R^2$, the function
$\psi_u(x) = (\psi(x),u)$ (an inner product computed with $\psi(x)$ treated as a vector in $\R^2$)
has no local maxima or minima at interior faces of $\Gd$. That is, every interior face $f$ (viewed
as an interior vertex in $\Gd'$) has neighbors $f_1$ and $f_2$ satisfying $\psi_u(f_1) \leq
\psi_u(f) \leq \psi_u(f_2)$. For generic $u$ (i.e., any $u$ whose slope is not parallel to one of
the $L_i$'s), the inequality can be made strict.

Now, to show that the $\{L_i\}$ form a T-graph, it remains only to show that they do not intersect
one another; while proving this, we will also show that $\psi(\Gd')$ partitions the convex polygon
$Q$ into convex polygons (the white faces). First, the maximal principle immediately implies that
$\psi(\Gd')$ lies in $Q$.  Furthermore, we claim that as one moves $x$ clockwise around each a
white interior face $w$ of $\Gd'$, $\psi(x)$ traces out a convex polygon in some fixed orientation
(clockwise or counterclockwise; we refer to this direction as the {\bf orientation of $w$} and
denote the polygon by $\psi(w)$). If this were not the case, then there would have to be vertices
$f_1,f_2,f_3,f_4$, in clockwise order around $w$ and some generic $u$ for which $\psi_u(f_1)$ and
$\psi_u(f_3)$ are less than both of $\psi_u(f_2)$ and $\psi_u(f_4)$. By the maximal principle, we
can find paths in $p_2$ and $p_4$ in $\Gd'$ from $f_2$ and $f_4$ to root vertices along which
$\psi_u$ is strictly increasing and paths $p_1$ and $p_3$ from $f_1$ and $f_3$ to root vertices
along which $\psi_u$ is strictly decreasing.   Now, let $p$ be a path in $\Gd'$ formed by
concatenating $p_1$ (reversed), a dotted line from $f_1$ to $f_3$, and $p_3$.  This path cannot
intersect $p_2$ or $p_4$ (since $\psi_u$ at any point on these two paths is greater than $\psi_u$
at any point on $p_1$ or $p_3$). However, the Jordan curve theorem implies that $p$ separates its
complement in $\Gd'$ into at least two connected components and that $f_2$ and $f_4$ (which lie on
either side of $p$ across the face $w$) are in separate components (this remains true even for the
graph $(\Gd^Q)'$ formed by adding to $\Gd'$ the edges connecting each cyclically consecutive pair
of outer vertices of $\Gd'$).  Now, the paths $p_2$ and $p_4$ both lead to root vertices at which
$\psi_u$ assumes a larger value than it does at any point along $p$, and these points are in the
same component of $(\Gd^Q)'$, a contradiction.  A similar argument shows that the outer faces $w$,
joined with this, have this orientation.  Another similar argument applies to black faces and shows
that as one moves $x$ around a black interior face $b$ of $\Gd'$, $\psi(x)$ traverses the
corresponding $L_i$ exactly once in each direction.

Next, we argue that all white faces have the same orientation.  It is enough to prove that any
white faces of $\Gd'$ (vertices of $\Gd$) $w_1$ and $w_2$ incident to a common black $b$ have the
same orientation.  Now, as $x$ traverses the boundary of the face $b$ in $\Gd'$, $\psi(x)$ traces
out the corresponding $L_i$ once in each direction; divide the faces incident to $b$ into two
categories according to the orientation of the edge shared with $b$. Clearly, if these faces do not
all have the same orientation, we can find two of them, $w_1$ and $w_2$ in opposite categories that
have opposite orientations.  In this case, $\psi(w_1)$ and $\psi(w_2)$ will lie on the same side of
$b$; let $u$ be vector orthogonal to $L_i$; assume without loss of generality that $\psi_u$ assumes
a larger value on points on $L_i$ than on other points of $w_1,w_2$. Let $f_1$ and $f_3$ be the
points in $\Gd'$ incident to $b$ whose images are the endpoints of $b$, and let $f_2$ and $f_4$ be
arbitrary points of $w_1$ and $w_2$ which do not lie on $b$. Let $p$ be formed by concatenating a
path $p_1$ from $f_1$ to a root on which $\psi_u$ is strictly increasing (reversed), a dotted line
from $f_1$ to $f_3$, and a path $p_3$ from $f_3$ to a root vertex along which $\psi_u$ is strictly
increasing; observing that $f_2$ and $f_4$ are on opposite sides of $p$, we derive a contradiction
through the Jordan curve argument described above.

Finally, suppose that two of the $L_i$ intersect.  Then there must be two faces $w_1$ and $w_2$ for
which $\psi(w_1)$ and $\psi(w_2)$ intersect.  The outer boundary of $(\Gd^Q)'$ is mapped with some
consistent orientation to $Q$.  Now, let $h:Q \rightarrow \Z$ at $x$ be the number of white faces
$\psi(w)$ which contain $x$ in their interiors.  It is clear that $h$ assumes the value $1$ near
the boundary.  We claim that $h$ is equal to one throughout $Q \backslash \psi (\Gd^Q)'$;
otherwise, there would be an $x$ in the interior of $Q$ (and not at the finitely many endpoints of
any $L_i$ or intersections of pairs of $L_i$) on the boundary of regions at which $h$ assumes
different values.  Such an $x$ must lie on some $L_i$, and it is not hard to see that the two white
faces incident to $x$ and $L_i$ must have opposite orientations.  \end{proof}

\subsection{Flat-face degeneracy} Now, suppose that $K$ is merely $1$-non-degenerate and not necessarily
$2$-non-degenerate; then we can formally construct $\psi$ exactly as above; in this case, however,
we cannot rule out that some of the $L_i$ may be parallel to one another---and in fact, some of the
$L_i$ may overlap. However, the same arguments given above still imply that for each white $w$,
$\psi(w)$ is {\it either} a convex face with some orientation (as described above) or a line
segment traversed once in each direction (like the black faces).  In the latter case, we say
$\psi(w)$ is a {\it degenerate face}.  In the presence of degenerate faces, we will consider
$\psi(w)$ and $\psi(b)$ to be incident to one another along an edge if and only if $w$ and $b$ are
adjacent vertices in $\Gd'$.

It is clear that if a white vertex $w$ is degenerate, then $\psi(b)$ is parallel to $\psi(w)$ for
each black $b$ adjacent to $w$.  A maximal component of the subgraph of $\Gd'$ consisting of
vertices on which $\psi$ is parallel to a given line is called a {\bf parallel component} of
$\Gd'$. Clearly,  the neighbor set of any white vertex in a parallel component is also in the
parallel component.

An {\bf extreme point} of a degenerate face $w$ is a vertex $f$ incident to $w$ for which $\psi(f)$
is an endpoint of $\psi(w)$.  The union of $\psi$-images of a parallel component is a segment which
we call an {\bf extended complete edge}. Now observe the following. \begin{lem} Each parallel
component $P$ is a $1$-cut. \end{lem}

\begin{proof} Observe that every $f$ which is an interior vertex of a black edge of $\Gd'$ in a
parallel cluster is the extreme vertex for the same number of black and white faces of $\Gd'$.  The
endpoints of the extended complete edge are extreme points of one more black vertex than white
vertices.  Since every face has exactly two extreme vertices, the result follows. \end{proof}

Similar arguments to those given in the proof of Theorem \ref{TgraphfromK} imply that as $x$
traverses the outside of a parallel component, $\psi(x)$ traverses the outside of the extended
complete edge exactly once in each direction.  Similar arguments to those of Theorem
\ref{TgraphfromK} imply that the extended complete edges form a T-graph. We say that $L = \{L_i \}$
forms a {\bf T-graph with overlaps} if $L_i$ satisfies all of the T-graph conditions except that
parallel pairs of $L_i$ are allowed to intersect (overlap) one another.   The above analysis
implies the following:

\begin{thm} \label{TgraphfromKdegenerate} Theorem \ref{TgraphfromK} still holds if $K$ is merely
$1$-non-degenerate and not necessarily $2$-non-degenerate---except that in this case, some of the
white faces may be degenerate (and so the T-graph may have overlaps).  Theorem \ref{corrthm} still
applies to T-graphs with overlaps.  \end{thm}

Even though some of the white faces are flat in the overlapping T-graph $\Gt$, we can define a dual
to the overlapping T-graph, containing these faces, using the graph structure of $\Gd$.  After
doing this, all of the arguments in the proof of Theorem \ref{corrthm} apply as before, so we still
have a martingale on the T-graph and have a measure preserving correspondence between spanning
forests and perfect matchings.

Recall that in any perfect matching, there is always exactly one
edge connecting a given $1$-cut to its complement, and that edge
contains a black vertex of the $1$-cut.  It is perhaps not
surprising that when we form the T-graph, $1$-cuts, in some sense,
play the same role as single black vertices. If we had simply
replaced all $1$-cuts in our original graph with single black
vertices, then, for a generic choice of weights, the T-graph would
not have any degenerate white faces.

\subsection{Extending the correspondence to degenerate weighted graphs}
\label{Tgraphmartingale}

Recall from Lemma \ref{degeneratelemma} that if we remove the
unused edges from $\Gd$, then the Kastelyn matrix for $\Gd$ is
$1$-non-degenerate (and hence Theorems \ref{corrthm}
\ref{TgraphfromKdegenerate} apply) for a generic choice of weight
functions $\nu$. Suppose, however, that the Kastelyn matrix for
$\Gd$ is not $1$-non-degenerate for a particular choice of weight
function $\nu$. Then we would like to take a generic sequence of
weights $\nu_i$ converging to $\nu$, look at the limit (or some
subsequential limit) of the corresponding T-graphs, and show that
the measure-preserving correspondence described in Theorem
\ref{corrthm} still holds for the limiting object. The problem is
that, as Figure \ref{nongeneric} makes clear, the limit of a
sequence of T-graphs need not be a T-graph at all; in fact, some
of the edge segments and faces may shrink to single points.

For practical computational applications, it may be sufficient to have the correspondence between
dimers and spanning forests for a generic choice of weights.  But a word of caution is in order.
Consider the dimer model whose T-graph is given by the right diagram in Figure \ref{nongeneric}; if
weights $\nu_i$ tend to a limit $\nu$ in such a way that the T-graphs have the graph on the left as
a limit, then the shrinking small triangle in the center of the diagram will become a ``trap'' for
the random walk on the T-graph, in that the expected amount of time that a walk spends on these
three vertices before exiting towards a root vertex tends to infinity; sampling algorithms that
rely on random walks will perform poorly for weights approximating $\nu$.  In this case, however,
one can simplify the limiting problem by reducing the three vertices in the small triangle at the
center to single vertex. The probability tends to one that only one of the ``long'' directed edges
(i.e., edges whose lengths are not tending to zero) extending outward from these three vertices
will appear in a random tree; given a spanning tree of the ``reduced'' graph, it is possible to
work out which ``short'' edges appear in the graph. The details of this and more general versions
of this reduction are left to the reader.

\section{Periodic and almost periodic T-graphs}\label{qp}

\subsection{Definitions for almost periodic T-graphs}

In this section, we prove some results about T-graphs which are motivated by the study of ergodic
Gibbs measures on tilings of infinite periodic
planar graphs.  More on this subject can be found in
\cite{KOS}, who cite the results of this section.
Our first aim here is to construct from periodic bipartite planar
graphs (and under certain conditions on the weights)
infinite T-graphs with a property called ``almost periodicity.''

Let $\Gd$ be embedded in the torus $\R^2 / \Z^2$ and let
$\Gd^\infty$ be the doubly periodic lift to $\R^2$ (we assume
that $\Gd^\infty$ is connected).  As before,
assume that $\Gd$ has $n$ white and $n$ black vertices.  Denote by
$v_{j,k}$ the vertex of $\Gd^{\infty}$ which lies in the square
$[j,j+1) \times [k,k+1)$ and whose projection to the torus is the
vertex $v \in \Gd$.  For the sake of simplicity
we will assume throughout this section that
$\Gd$ has no unused edges and that it has generic weights.
The non-generic weight case requires a slightly finer analysis
which we choose not to go into here. See however \cite{KOS}.

A function $f$ on the vertices of $\Gd$ is {\bf $(\alpha, \beta)$-periodic} if $f(v_{j+x,k+y}) =
\alpha^x\beta^y f(v_{j,k})$ for all $(v_{j,k}) \in \Gd^{\infty}$. Say $f$ is {\bf almost periodic}
if it is $(\alpha,\beta)$-periodic and $\alpha$ and $\beta$ have modulus one (but are not
necessarily roots of unity). In this case, we write $\alpha = e^{2 \pi i a}$ and $\beta = e^{2 \pi
i b}$.  If $a$ and $b$ are rational, then $f$ is doubly periodic with some period.

For a fixed $(\alpha,\beta)$ the linear space of $(\alpha,\beta)$-periodic functions is
$2n$-dimensional and is parametrized by the space of functions on one period of
$\Gd^{\infty}$---which we can represent as a single copy of $\Gd$.
It has a natural basis consisting of functions $\delta_v$ whose value is $1$ at $v\in[0,1)^2$
and zero at other vertices in the fundamental domain.
Let $K$ be a Kastelyn matrix for
$\Gd$ and $K^{\infty}$ an infinite-dimensional Kasteleyn matrix for $\Gd^{\infty}$ which is a lift
of $K$.  We can think of $K^{\infty}$ as a linear function from the set of functions on the black
vertices of $\Gd^{\infty}$ to functions on the white vertices of $\Gd^{\infty}$.  Since this
function maps $(\alpha,\beta)$-periodic functions to $(\alpha,\beta)$-periodic functions, it
induces a linear map from the $n$-dimensional space of functions on the black vertices of $\Gd$ to
the $n$-dimensional space of functions on white vertices of $\Gd$.
Denote by $K_{\alpha,\beta}$ the matrix of this linear map in the basis $\{\delta_v\}$.

\begin{figure} \begin{center} \leavevmode \epsfbox[86 535 276 717]{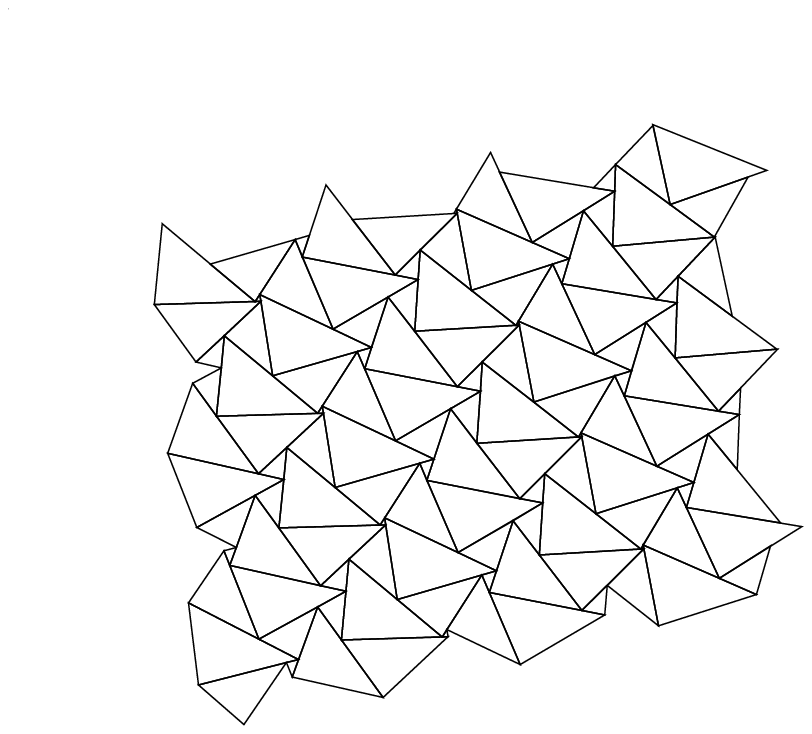} \caption{\label{righttri}
Almost periodic T-graph mapping of the honeycomb graph, with periodic edge weights $4,5,$ and $6$
according to direction.  The edges of the graph shown correspond to black vertices of the honeycomb
lattice; the triangular faces of the graph shown correspond to white vertices of the honeycomb
lattice.} \end{center} \end{figure}

The determinant $\det K_{\alpha,\beta}$ is a polynomial function of $\alpha$ and $\beta$; in
particular for certain $(\alpha,\beta)$ (corresponding to zeros of this polynomial function) the
matrix $K_{\alpha,\beta}$ has a non-trivial null space, and hence we can find
$(\alpha,\beta)$-periodic functions $f$ and $g$ satisfying $K^\infty f=0$ and $gK^\infty=0$.  If
the polynomial $\det K_{\alpha, \beta}$ {\it happens} to have a zero $(\alpha, \beta)$ that lies on
the unit torus of complex variable pairs that both have modulus one, then $f$ and $g$ are almost
periodic.  If, furthermore, $f$ and $g$ {\it happen} to be nowhere zero, then we can define an
infinite T-graph as follows.  First, observe that the function $\tilde{K}^{\infty}_1(vw) =
f(v)g(w)K^{\infty}(v,w)$ on edges $vw$ of $\Gd^\infty$ is a nowhere zero flow. The dual of this
flow is the gradient of a function $\psi_1$ on $\Gd'$. Similarly the dual of
$\tilde{K}^{\infty}_2(vw) =  f(v) \overline{g(w)}K^{\infty}(v,w)$ is the gradient of a function
$\psi_2$ on $\Gd'$ (where $\overline{g}$ denotes the complex conjugate of $g$).  We may assume
(multiplying $g(w)$ by a generic modulus one complex number if necessary) that $g(w) +
\overline{g(w)} = 2\Re g(w)$ is also nowhere zero.  Then we can think of $\tilde{K} = \tilde{K}_1 +
\tilde{K}_2$ as an infinite Kastelyn matrix and $\psi=\psi_1+\psi_2$ as the corresponding T-graph.
We will call a mapping $\psi$ from $(\Gd^\infty)'$ to $\R^2$, constructed in this way, an {\bf
almost periodic T-graph mapping}. See Figure \ref{righttri}.

We remark that, given a fixed $\nu$, the number of $\alpha,\beta$ on the unit torus for which $\det
K_{\alpha,\beta} = 0$ also plays a fundamental role in \cite{KOS}, where it is shown that the
minimal specific free energy ergodic Gibbs measure on perfect matchings of the infinite weighted
graph $\Gd^{\infty}$ is {\it smooth} if the corresponding polynomial $K_{\alpha, \beta}$ polynomial
has $0$ roots on the unit torus and {\it rough} if it has $2$ roots (necessarily complex
conjugates) on the unit torus (in the non-generic case of a single root, it is rough only when
$\frac{d}{d\alpha}\det K_{\alpha,\beta}=\frac{d}{d\beta}\det K_{\alpha,\beta}=0$).
The terms ``smooth'' and ``rough'' come from the statistical
physics literature and are defined in \cite{KOS}.  The main goal of this section is to
prove that when the choice of
weights is generic, the number of modulus-one values of $(\alpha, \beta)$ that are roots of $\det
K_{\alpha,\beta}$ always belongs to the set $\{0,2\}$.

\subsection{Generic points on the variety of almost periodic T-graphs}

Write $\mathbb R_+$ for the set of strictly positive real numbers, $\mathbb C_+$ for the set of
non-zero complex numbers, and write $\mathbb P^k$ for $k$-dimensional complex projective space.
Suppose that $|V| = 2n$, and define a variety $X \subset \mathbb R_+^{|E|} \times \mathbb C_+^2
\times \mathbb P^{n-1}$ by:
$$X = \{(\nu, \alpha, \beta, f) : K_{\alpha, \beta} f = 0 \}$$
Here $f$ is an element in $\mathbb P^{n-1}$, which is a one-dimensional subspace of $\mathbb C^n$,
and by $K_{\alpha, \beta} f = 0$ we mean that this subspace lies in the null space of $K_{\alpha,
\beta}$.  By abuse of notation, if $f$ is a non-zero function on the black vertices of $\Gd$, we
will also use $f$ to denote the element of $\mathbb P^{n-1}$ given by the linear span of $f$.
Denote by $\tilde X$ the subset of $X$ consisting of points for which $|\alpha| = |\beta| = 1$.
Denote by $\Adj(K_{\alpha, \beta})$ the {\bf adjugate} matrix of $K_{\alpha, \beta}$, whose entries are
the $(n-1) \times (n-1)$ minors of $K_{\alpha,\beta}$ (so that $K_{\alpha,\beta} \Adj(K_{\alpha,
\beta}) = \det K_{\alpha,\beta}$).  It is easily seen that $\Adj(K_{\alpha,\beta})$ is identically equal
to zero if and only if the rank of $K_{\alpha,\beta}$ is less than $n-1$; and if the rank of
$\Adj(K_{\alpha,\beta})$ is exactly $n-1$, then at least one column of $\Adj(K_{\alpha,\beta})$ is a
non-zero vector whose span is the null space of $\Adj(K_{\alpha,\beta})$.  The following is the main
result of this section:

\begin{thm} \label{generic2to1} The variety $\tilde X$ is irreducible.  For a generic choice of
$\nu$, there are either zero or two quadruples $(\nu, \alpha, \beta, f)$ in $X$.  When the latter
is the case and $\nu$ is generic, then the corresponding $\alpha,\beta,f$ are such that
$\Adj(K_{\alpha,\beta})$ has rank $n-1$, all of its coordinates are non-zero, and $f$ is
given by any column of the (rank one) matrix $\Adj(K_{\alpha,\beta})$. \end{thm}

Let us say a word about the significance of this theorem to T-graph classification before we prove
it.  By obvious symmetry, Theorem \ref{generic2to1} implies that that for generic $\nu$, there are
either zero or exactly two quintuples $(\nu, \alpha, \beta, f, g)$ with $gK_{\alpha,\beta} = 0$ and
$fK_{\alpha,\beta}$. Recall that our almost periodic T-graphs were defined to have gradient given
by $\tilde{K}^{\infty}(vw) = 2f(v)\Re(g(w))K^{\infty}(v,w)$. Since $f$ and $g$ are uniquely
determined up to complex conjugacy and multiplication by a constant factor, this implies that the
almost periodic T-graph is completely determined up to rotations (which arise from multiplying $f$
by a modulus one constant), constant rescalings (which arise from multiplying either $g$ or $f$ by
a real constant), reflection (which comes from complex conjugacy), translations of the image space
(which arise from the fact that $\tilde{K}^\infty$ only determines the T-graph mapping up to an
additive constant) and ``translation of the domain, or a limit of such translations.'' To explain
the last symmetry, note that multiplying both $f$ and $g$ by $\alpha^m\beta^n$ is equivalent to
composing the T-graph mapping with translation of the domain by $(m,n)$.  If one of $\alpha,\beta$
is irrational, then we can achieve any modulus one number as a limit of numbers of the form
$\alpha^m \beta^n$.  We summarize these observations informally by saying that ``the almost
periodic T-graph mapping corresponding to $\nu$ is unique up to affine orthogonal transformations
of the image and translations of the domain.'' We say two T-graphs are equivalent if one can be
obtained from the other via a symmetry of this sort. Note, of course, that if $\alpha$ and $\beta$
are both rational, then multiplying $f$ and $g$ by a modulus one number is {\em not} necessarily
the same as a domain translation, or even a limit of such translations.  In this case, there is a
one parameter family of T-graph equivalency classes.

We will now prove Theorem \ref{generic2to1} in stages, beginning with the following lemma.  First,
denote by $X'$ the projection of $X$ onto its first three coordinates $(\nu,\alpha,\beta)$; i.e.,
$X'$ is the zero set of the polynomial $P(\nu,\alpha,\beta) = \det K_{\alpha, \beta}$.

\begin{lem} \label{genericnozeroentry} The variety $X'$ is irreducible.  Moreover, for a generic
point $(\nu,\alpha,\beta)$ on $X'$, the matrix $\Adj(K_{\alpha,\beta})$ has no zero entries, and the
$f$ for which $(\nu,\alpha,\beta,f) \in X$ is unique. \end{lem}

\begin{proof} Clearly, $P$ is affine linear as a function of $\nu(e)$, that is $P=\nu(e) P_e
+P_e',$, where $P_e$ and $P_e'$ do not involve $\nu(e)$. If we could write $P=P_1P_2$, then each
$\nu(e)$ must occur in either $P_1$ or $P_2$, but not both.  Since the multiplicity of the $\nu(e)$
terms determine the multiplicity of $\alpha$ and $\beta$ in each monomial, this implies that there
is no cancellation when multiplying out $P_1$ times $P_2$ (i.e., there are no monomials that can
represented as a product of a monomial in $P_1$ and a monomial in $P_2$ in two different ways).
Thus, each monomial in $P_1$ times a monomial of $P_2$ corresponds to a matching. Let $E_1,E_2$ be
the set of edges represented in $P_1,P_2$, respectively, and $V_1,V_2$ their vertices. If an edge
$e$ connected a vertex $v_1$ of $V_1$ to a vertex $v_2$ of $V_2$, then its weight could not occur
in either $P_1$ or $P_2$, since if it occurred in a monomial of, say, $P_1$, then the product of
that monomial with a monomial of $P_2$ that included a factor of $\nu(e')$ with $e'$ incident to
$v_2$ (such a monomial exists by definition) would {\em not} correspond to a matching, since it
would involve two edges incident to $v_2$. Thus $e$ must be unused, a contradiction. Thus, if $P =
P_1 P_2$, then one of the $P_i$---say, $P_2$---must be a function of $\alpha$ and $\beta$ alone.
Since each combination of edge weights corresponding to a matching occurs in exactly one monomial
of $P$, we conclude that $P_2$ is a monomial in $\alpha$ and $\beta$.

Furthermore $P$ is irreducible when considered as a polynomial in both the edge weights and
$\alpha,\beta$, except for a monomial factor in $\alpha$ and $\beta$. That is, if
$P=P_1(\nu,\alpha,\beta) P_2(\nu, \alpha,\beta)$ then one of the $P_i$ consists of a single
monomial in $\alpha$ and $\beta$. To see this, note that by the previous result, we may assume
without loss of generality that $P_2$ is a polynomial in $\alpha, \beta$ alone; and since we are
assuming $\alpha \not = 0, \beta \not = 0$, the variety is not changed if we divide out by this
term so that $P$ is an irreducible polynomial.

Fix an edge $e$ and consider the polynomial $P_e$ as defined above. Since $P$ is irreducible and
$e$ occurs in a proper subset of the set of all matchings, the zero set of $\nu(e) P(e)$,
intersected with $X'$, forms a proper subvariety of $X'$.  In other words, on a generic
subset of $X'$, none of the entries of $\Adj(K_{\alpha,\beta})$ corresponding to an edge in $\Gd$ are
zero.  Since $\Adj(K_{\alpha,\beta})$ has rank at most one, and every row and column has a non-zero
entry, we conclude that every entry of $\Adj(K_{\alpha, \beta})$ is non-zero and $f$ is the span of any
column of $\Adj(K_{\alpha,\beta})$. \end{proof}

 \begin{lem} \label{genericeverynozeroentry} For a generic choice of weights $\nu$,
{\it every} pair $\alpha, \beta$ for which $(\nu, \alpha, \beta) \in X$ is such that
$\Adj(K_{\alpha,\beta})$ has no zero entries, and the $f$ for which $(\nu,\alpha,\beta,f) \in X$ is
unique. \end{lem}

\begin{proof} Lemma \ref{genericnozeroentry} implies that for {\it generic} edge weights $\nu$, $P$
and $P_e$ have no common factor as functions of $\alpha$ and $\beta$ except for monomial factors.
To see this, by irreducibility note that there exist polynomials $Q_1=Q_1(\alpha,\beta,w)$ and
$Q_2=Q_2(\alpha,\beta,w)$ such that $PQ_1+P_eQ_2=Q(\alpha,w)$ where $Q$ is a nonzero polynomial
depending only on $\alpha$ and the weights $w$, not on $\beta$. Similarly there exist $Q_3,Q_4$
such that $PQ_3+P_eQ_4=Q'(\beta,w)$ where $Q'$ is a nonzero polynomial independent of $\alpha$.
Plugging in generic values for $w$, $Q$ and $Q'$ will still be nonzero, but any common factor of
$P$ and $P_e$ is a common factor of $Q$ and $Q'$ which is impossible. So $P$ and $P_e$ have no
common factor for generic $w$.

Therefore, when $\nu$ is fixed generically, by Bezout's theorem $P$ and $P_e$---viewed as
polynomials in $\alpha$ and $\beta$---have a finite number of common zeros. By genericity none of
these zeros lies on the unit torus (since for any positive real $x$, we can choose $\nu_x$ so that
$P(\nu,\alpha,\beta) = P_{\nu_x, x\alpha,x\beta}$; and replacing $\nu$ with such a $\nu_x$, for a
generic choice of $x$, preserves the genericity of the weights). \end{proof}


\begin{lem} \label{unbounded} Any almost periodic T-graph
mapping $\psi$ is unbounded as a function of $(\Gd^{\infty})'$.
Moreover if $u$ is any vector in $\R^2\backslash\{0\}$, then
$(\psi,u)$ is unbounded if it is not identically equal to a
constant (in which case $\psi$ is degenerate---i.e., its image is
contained in a line).
\end{lem}

\begin{proof} Suppose that $f$ is $(\alpha,\beta)$-periodic and
$g$ is $(\gamma,\delta)$-periodic with $\gamma = e^{2\pi i c}$ and
$\delta = e^{2 \pi i d}$. Then $\tilde K^{\infty}_1 (v_{j,k},
w_{j+\ell,k+m})$ is a function of $\ell,m$ whose real and
imaginary parts can both be written in the form $\cos(a\ell + bm +
x)\cos(c\ell + dm + y)$ times a constant, for some $x$ and $y$.

If $\psi$ were bounded on $\Gd'$, then the corresponding
martingale on the T-graph would almost surely converge (by the
martingale convergence theorem), and there would thus have to be a
path of vertices $v_1,v_2, \ldots$ for which $\psi(v_i)$ converges
to a constant. We claim that this is impossible. It is enough to
show that for some $\epsilon$, the set of edges $(vw)^*$ for which
$0<\tilde K^\infty(v,w) < \epsilon$ has no infinite cluster.  For
some $N > |\Gd|$, we can always find $\epsilon$ small enough so
that the distance between any two clusters of $(\ell,m) \in \Z^2$
(viewed as points in $\Z^2$) on which $0<\cos (a\ell + bm + x) <
\epsilon^{1/2}$ is at least $2N$ times the diameter of the largest
such cluster, and similarly for clusters on which $0<\cos (c\ell +
dm + y) < \epsilon^{1/2}$.  (This is trivial if $a$ and $b$ are
rational, since the function is periodic in that case; if they are
irrational, then we can find $\epsilon_0$ for which there is no
integer pair $(n_1,n_2)$ for which $n_1a + n_2b$ is less than
$\epsilon_0$ (modulo $2\pi$) and $|n_1+n_2| \leq 2N$. Choose
$\epsilon_0$ small enough that there can't be two values differing
by $\epsilon_0$ (modulo $2\pi$) with cosines $\epsilon$ apart.)
Now, it is clear that the largest cluster of $\ell,m$ on which
even one of these statements holds is at most $2N$; since the
gradient of $\psi$ has norm at least $\epsilon$ when neither
statement holds, we conclude that $\psi$ cannot be bounded.

The same argument shows that there cannot exist a non-zero vector
$u \in \R^2$ for which the inner product $(\psi(v),u)$ is bounded
as a function of $v$, unless $(\psi(v),u)$ is constant.
\end{proof}

\begin{lem} \label{atmosttwo} If $\nu$ is generic, then the maximum number of linearly independent,
almost periodic solutions to $K^{\infty}f = 0$ (or similarly, solutions to $gK^{\infty} = 0$) is
two. If there are two solutions, which are $(\alpha,\beta)$- and $(\gamma,\delta)$-periodic, then
$\alpha=\bar\gamma$ and $\beta=\bar\delta$. \end{lem}

\begin{proof} For each $\alpha$ and $\beta$, the left null space
of $K_{\alpha,\beta}$ has the same dimension as the right null
space.  Now, suppose that $f$ is $\alpha,\beta$-periodic and $g$
is $\gamma,\delta$-periodic with $\gamma = e^{2\pi i c}$ and
$\delta = e^{2 \pi i d}$.  Then as in the proof of Lemma
\ref{unbounded}, $\tilde K^{\infty} (v_{j,k}, w_{j+\ell,k+m})$ is
a function of $\ell,m$ whose real and imaginary parts can both be
written in the form $\cos (a\ell + bm + x)\cos (c\ell + dm + y)$
times a constant, for some $x$ and $y$.  Let $S$ be a cycle in
$\Gd'$; if we lift it to $(\Gd^{\infty})'$, then its endpoints are
its starting points plus an integer pair, $(n_1,n_2)$.  Now, we
would like to determine the asymptotics of $\psi_1$ and $\psi_2$
(whose derivative is the dual of $\tilde K^{\infty}$) along
$S^{\infty}$ (a periodic lifting of $S$ to $\Gd^{\infty}$).
Expanding the cosines in exponentials, this involves adding up
$|S|$ separate sequences (functions of $\ell$) of the form:
$$\sum_{\ell=1}^{n_1} e^{2 \pi i [(x+\ell a) \pm (y+\ell c)]}$$
and $|S|$ sequences of the corresponding form for $m$.

Clearly, $\psi$ will remain bounded independently of $x$ and $y$,
provided $a \neq \pm c\bmod 2 \pi$ and $b\neq\pm d\bmod 2\pi$. In
fact we must take the same sign for both equalities: unless $(a,b)
= \pm (c,d)\bmod 2 \pi$ it is possible to find an independent pair
of integer vectors $(m_1,n_1)$ and $(m_2,n_2)$ for which $am_1 +
bn_1 \not = \pm(cm_1 + dn_1)\bmod 2\pi$ and similarly $am_2 + bn_2
\not = \pm(cm_2 + dm_2)\bmod 2 \pi$. Taking $S_1$ and $S_2$ to be
corresponding paths, we may deduce that $\psi$ is bounded unless
$(\alpha,\beta)$ and $(\gamma,\delta)$ are either equal to one
another or conjugates; by Lemma \ref{unbounded} $(\alpha,\beta)$
and $(\gamma,\delta)$ are either equal to one another or
conjugates.

Now suppose we have $(a,b) = \pm (c,d)$.  Then for the sums
corresponding to steps in $S$,
$$\sum_{\ell=1}^{n_1} \cos (x + \ell a) \cos (y \pm \ell a)$$ is
{\bf approximately linear} as a function of $n_1$, that is, equal
to a linear function plus a bounded function.  If there were three
linearly independent solutions $f_1,f_2,f_3$ to $Kf = 0$, and
$\psi_1, \psi_2,\psi_3$ are formed using $g$ and $f_1,f_2,f_3$,
then a linear combination of the $\psi_1,\psi_2,\psi_3$ would be
approximately the linear function zero (i.e., bounded), a
contradiction, by Lemma \ref{unbounded}.

Finally, since it is clear that $(\alpha,\beta)$ is not real (i.e., not equal to $\pm 1$) for a
generic choice of $\nu$, so any almost periodic $f$ or $g$ will be a strictly non-real function, that
is, linearly independent from its complex conjugate, which is also a zero of $K^{\infty}$.
\end{proof}

Now, Theorem \ref{generic2to1} now follows immediately from Lemma \ref{genericeverynozeroentry} and
Lemma \ref{atmosttwo}.

\old{
\subsection{Degenerate almost periodic T-graphs and non-generic points of variety}

We have described the almost periodic T-graphs for a generic choice of weights; but the non-generic
cases are also of interest.  Of particular interest are the choices of weights $\nu$ for which
$\det K_{\alpha, \beta} = 0$ for $\alpha, \beta \in \pm 1$, since these give rise to periodic (not
merely almost periodic) T-graphs.

The proofs of the almost periodic T-graph classification in the previous section excluded three
types of non-generic points on $X$:

\begin{enumerate}
\item Points $(\nu,\alpha,\beta, f)$ for which $\alpha$ and $\beta$ are real.
\item Points $(\nu,\alpha,\beta,f)$ for which $K_{\alpha,\beta}^*$ has at least one zero entry but is not identically zero.
\item Points $(\nu, \alpha, \beta, f)$ for which {\it all} entries of $K_{\alpha,\beta}$ are zero (i.e.,
the rank of $K_{\alpha,\beta}$ is less than $n-1$). \end{enumerate}

Our first observation is that these particular ``non-generic'' points may all be obtained as limits
of the generic ones.

\begin{lem} The points in $(\nu,\alpha,\beta,f) \in X$ at which $\alpha, \beta \not \in \mathbb R$
and $K_{\alpha,\beta}^*$ has all non-zero entries form a dense subset of $X$ in the Euclidean
topology. \end{lem}

\begin{proof} \end{proof}

This fact enables us to conclude almost periodic T-graph uniqueness in the case of the first two
symmetry types:

\begin{thm} If $(\nu,\alpha,\beta,f) \in X$ is such that $K_{\alpha,\beta}$ has rank $n-1$, then
there are exactly two points in $X$ first first coordinate $\nu$ (complex conjugates of one
another) if at least one of $\alpha,\beta$ is not real, and exactly one such point if
$\alpha,\beta$ are real.  In the latter case, the image of the corresponding T-graph mapping is
completely contained in a line. \end{thm}

The latter symmetry is more special.  First, it can only appear when $\alpha$ and $\beta$ are both
real:

\begin{thm} Whenever $(\nu,\alpha,\beta,f) \in X$ is such that the rank of $K_{\alpha,\beta}$ is
less than $n-1$, and $\alpha,\beta$ have modulus one, then we have \begin{enumerate} \item Both
$\alpha$, and $\beta$ are real, i.e., equal to $\pm 1$.
\item The rank of $K_{\alpha,\beta}$ is exactly $n-2$.
\item The set of points in $X$ with first coordinate given by $\nu$ is such that $\alpha$ and $\beta$
are completely determined, but $f$ may be any line contained in the two-dimensional null space of
$K_{\alpha,\beta}$.  \end{enumerate}\end{thm}

The following lemmas aid us in proving this:

\begin{lem} The $T$-graph corresponding to $\nu, \alpha, \beta, f, g$ is not identically zero for
any choice of $f$ and $g$; for every choice of $f$ and $g$ it is approximately linear and
unbounded. \end{lem}

\begin{lem} We have $\nabla P(\nu, \alpha, \beta) = 0$ (where the gradient is of $P$ when $P$ is
viewed as a function of $\alpha$ and $\beta$ alone) if and only if the rank of $K_{\alpha, \beta}$
is $n-2$. If the gradient is non-zero, then the rank of $K_{\alpha, \beta}$ is $n-1$. \end{lem}
}

\old{
\subsection{Periodic T-graphs and drift}

Now, suppose that $\alpha$ and $\beta$ are real (i.e. in $\{\pm1\}$ and $\det K_{\alpha,\beta} = 0$
and has nowhere zero vectors in its left and right null spaces (so T-graphs are well-defined).  It
is natural to wonder whether there is any difference between cases when $K_{\alpha,\beta}$ has rank
$n-1$ and when it has rank $n-2$. This turns out to be related to whether the random spanning tree
model which corresponds to $K^{\infty}$ has drift. More precisely, this means that in the
stationary distribution of the random walk on the torus, the expected horizontal or vertical change
in a step of the walk is non-zero.

\begin{thm} If $\alpha, \beta\in\R$ and $\det K_{\alpha,\beta} = 0$, then there is a fully
non-degenerate T-graph (i.e., a T-graph mapping approximating an invertible linear function) if and
only if the rank of $K_{\alpha,\beta}$ is $n-2$. If the rank is $n-1$, then there exists a T-graph
which is contained in a single line.  The random walk on the periodic weighted graph
$(\Gd^{\infty})'$ has a net drift if and only if the rank is $n-1$. \end{thm}

\begin{proof} Suppose that there exists only one solution
to $\det K_{\alpha,\beta} = 0$, and so the T-graph is degenerate,
but suppose that there is no net drift. For $x\in(\Gd^{\infty})'$
let $\phi_0(x)$ be the limit in $N$ of the expected value of the
position of the $N$-step random walk (viewed as walk on the $\R^2$
positions of vertices of $\Gd'$ as it is embedded doubly
periodically in $\R^2$). Since there is no net drift, this limit
exists for every vertex of $(\Gd^{\infty})'$, and $\phi_0$ is a
non-degenerate approximately linear function. However, $\phi_0$ is
harmonic and so also defines a non-trivial two-dimensional,
periodic T-graph. This is a contradiction.
\end{proof}

Using the language and results of \cite{KOS}, this states that the
difference between so-called smooth and rough phases of ergodic
dimer Gibbs measures corresponds to the difference between uniform
spanning tree models on lattice graphs with and without drift. The
appearance of smooth phases is also related in \cite{KOS} to the
appearance of crystal facets in random surface models based on
height functions of perfect matchings.  This surprising connection
between crystal facet appearance and T-graph drift is intriguing
and merits further study.
}

\end{document}